\documentclass[bibay1]{asl}

\usepackage{mathrsfs}
\usepackage{braket}

\newtheorem{thm}{Theorem}
\newtheorem{lem}[thm]{Lemma}
\newtheorem{cor}[thm]{Corollary}
\theoremstyle{definition}
\newtheorem{defin}[thm]{Definition}

\newcommand{\scott}[1]{\mathscr{#1}}
\newcommand{\PA}{\mathrm{PA}}
\newcommand{\WKL}{\mathrm{WKL}_0}
\newcommand{\ACA}{\mathrm{ACA}_0}
\newcommand{\RCA}{\mathrm{RCA}_0}

\newcommand{\Con}{\mathop\mathrm{Con}}

\newcommand{\SSy}{\mathop\mathrm{SSy}}
\newcommand{\cod}{\mathop\mathrm{set}\nolimits}
\newcommand{\Cod}{\mathop\mathrm{Cod}\nolimits}
\newcommand{\rep}{\mathop\mathrm{rep}}
\newcommand{\Th}{\mathop\mathrm{Th}}
\newcommand{\La}{\mathscr{L}}

\newcommand{\prf}{\vdash}
\newcommand\nprf\nvdash

\newcommand{\embin}{\prec}
\renewcommand{\set}[1]{\Set{#1}}
\newcommand{\power}{\mathscr{P}}

\title{A note on standard systems and ultrafilters}
\author{Fredrik Engstr\"om}
\revauthor{Engstr\"om, Fredrik}
\address{Department of Philosophy \\ G\"oteborg University \\ 
  Box 200, 405 30 G\"oteborg, Sweden}
\email{fredrik.engstrom@filosofi.gu.se}

\begin{document}
\begin{abstract}
Let $(M,\scott X) \models \ACA$ be such that $P_\scott X$, the collection of all
unbounded sets in $\scott X$, admits a definable complete ultrafilter and let
$T$ be a theory extending first order arithmetic coded in $\scott X$ such that
$M$ thinks $T$ is consistent. We prove that there is an end-extension $N \models
T$ of $M$ such that the subsets of $M$ coded in $N$ are precisely those in
$\scott X$. As a special case we get that any Scott set with a definable
ultrafilter coding a consistent theory $T$ extending first order arithmetic is
the
standard system of a recursively saturated model of $T$.
\end{abstract}

\maketitle

The standard system of a model $M$ of PA (the first order formulation of Peano
arithmetic) is the collection of standard parts of the parameter definable
subsets of $M$, i.e., sets of the form
$X \cap \omega$, where $X$ is a parameter definable set of $M$, and
$\omega$ is the set of natural numbers. It turns out that the standard
system tells you a lot about the model; for example, any two countable
recursively saturated models of the same completion of PA with the same
standard system are isomorphic. A natural question to ask is then which
collections of subsets of the natural numbers are standard systems. This problem
has become known as the \emph{Scott set problem}.
 
For countable models the standard systems are exactly the countable Scott sets,
i.e., countable boolean algebras of sets of natural numbers closed under
relative recursion and a weak form of K\"onig's lemma \cite{Scott:62}.
It follows, by a
union of chains argument (see \cite{Knight.Nadel:82}), that for models of
cardinality at most $\aleph_1$ the
standard systems are exactly the Scott sets of cardinality at most $\aleph_1$.

If the continuum hypothesis holds this settles the Scott set problem. However,
if CH fails then very little is known about standard systems of models of
cardinality strictly greater than $\aleph_1$, although it is easy to see that
any standard system of any model is a Scott set (see \cite{Kaye:91}).

The key notion of a definable ultrafilter plays a central role in
this paper. Suppose $\scott{X}$ is a collection of
sets of natural numbers, and $P_\scott{X}$ is
the collection of infinite members of $\scott{X}$. A filter
$F$ on $P_\scott{X}$ is said to be definable if for all $A \in \scott X$
$\set{a  \in \omega  | (A)_a \in F} \in \scott X$,  where $(A)_a=\set{b \in
\omega | \langle a,b \rangle \in A}$ and $\langle \cdot,\cdot \rangle$ is some
canonical pairing function. The central theorem of this paper is Theorem
\ref{maintheorem} below (whose proof will be presented in section
\ref{sec:con}).

\begin{thm}\label{maintheorem}
Let $\scott X$ be a Scott set that carries a definable ultrafilter and let $T
\in
\scott X$ be a consistent completion of PA. Then there is a recursively
saturated model of $T$ with standard system $\scott X$.
\end{thm}

It should be noted that all Scott sets that carry a definable ultrafilter are
arithmetically closed (see Lemma \ref{lem:def.ar}) and that any countable
arithmetically closed Scott set carries a definable ultrafilter, thus our
result partly extends Scott's theorem (Theorem \ref{thm:scott} below).

Our hope is that this paper may shed some light on  the Scott set problem for
uncountable standard systems. We will end the paper by offering some
hints as to how this might be done.

Some of the results in this paper can also be found in the author's doctoral
thesis \cite{Engstrom:04}.

\section{Background and preliminaries}

Let $\La_A$ be the language of arithmetic with symbols $0,1,+,\cdot,<$ and
$\PA$ the ordinary axiomatization of arithmetic in first order logic. Let
$\mathbb N$ be the standard model of PA and $\omega$ be its domain, i.e., the
set of natural numbers. If $a\in M\models\PA$ let $I^M_{\mathord < a}$ (or
$I_{\mathord < a}$ if $M$ is understood from the context) be the initial segment
$\set{b \in M | M \models b < a}$ of $M$. $I \subseteq_e M$ means that $I$ is an
initial segment of $M$, i.e., $I^M_{\mathord < a}\subseteq I$ for all $a \in I
$. $I \subseteq_e M$ is a cut if it is closed under taking successors.
Let $x \in y$ be a $\Delta_0$-formula in the language of arithmetic used for
coding sets in $\mathrm{I}\Sigma_1$ (see
\cite{Kaye:91}).\footnote{$\mathrm{I}\Sigma_1$ is PA with
the induction axioms restricted to $\Sigma_1$ formulas.}
If $I$ is a cut in $M$
and $a \in M$ then $\cod_{M/I}(a) = \set{b \in I | M \models b \in a}$. The set
of all coded subsets of $I$ will be denoted $\Cod(M/I)= \set{\cod_{M/I}(a) | a
\in M}$. The standard system of $M$ is $\SSy(M)=\Cod(M/\omega)$.

If $I \subseteq M$ is a cut closed under multiplication and  $\scott X \subseteq
\power(I)$ we can regard the structure $(I,\scott X)$ as a second order
structure in the language $\La_A$. The subtheories $\RCA$, $\WKL$ and $\ACA$ of
second order arithmetic are as defined in \cite{Simpson:99}. A cut $I \subseteq
M$ is semiregular iff $(I,\Cod(M/I)) \models \WKL$, and strong iff
$(I,\Cod(M/I)) \models \ACA$. 

When we write $(M,\scott X)$ it will be understood that $(M,\scott X)$ is a
second order structure in the language $\La_A$, i.e., that $M$ is equipped with
addition and multiplication and $\scott X \subseteq \power(M)$. $(M,\scott X)$
is an $\omega$-model if $M=\mathbb N$ is the standard model of PA. In this case
$(\mathbb N,\scott X) \models \WKL$ iff $\scott X$ is a Scott set. We will often
specify an $\omega$-model by only specifying $\scott X \subseteq
\power(\omega)$.

Given a theory $T$, in the language of arithmetic, we say that a set $A
\subseteq \omega$ is represented in $T$ if there is a formula
$\varphi(x)$ such that  $T \prf \varphi(n)$ for all $n \in A$ and $T \prf \lnot
\varphi(n)$ for all $n \in \omega \setminus A$. By $\rep(T)$ we denote the
collection of all sets represented in $T$. In \cite{Scott:62} Scott proved a
variant of the following theorem:

\begin{thm}\label{thm:scott}
For any  countable $\omega$-model $\scott{X}$ and complete consistent theory $T
\supseteq \PA$ the following are equivalent:
\begin{itemize}
\item[$(i)$] $\scott{X} \models \WKL$ and $\rep(T) \subseteq \scott X$, and
\item[$(ii)$] there is a countable nonstandard model $M \models T$ with
$\SSy(M)=\scott{X}$. 
\end{itemize}
\end{thm}

This theorem could be taken a bit further: Given a Scott set $\scott{X}$ of
cardinality $\aleph_1$ there is a nonstandard model of $\PA$ with standard
system $\scott{X}$: \footnote{In \cite{Smorynski:84}, Smor\'ynski gives
Ehrenfeucht and Jensen \cite{Ehrenfeucht.Jensen:76} and independently Guaspari
\cite{Guaspari:79} the honor of the main lemma used to prove the theorem. He
also writes that the observation that the theorem follows from the lemma is due
to Guaspari in \cite{Guaspari:79}. It is not clear to us if this is correct.
What is clear is that the lemma and the theorem appears explicitly in Knight and
Nadel \cite{Knight.Nadel:82}.}

\begin{thm}\label{union}
For any  $\omega$-model $\scott{X}$ of cardinality at most $\aleph_1$ and
complete consistent theory $T \supseteq \PA$ the following are equivalent:
\begin{itemize}
\item[$(i)$] $\scott{X} \models \WKL$ and $\rep(T) \subseteq \scott X$, and
\item[$(ii)$] there is a nonstandard model $M \models T$ with
$\SSy(M)=\scott{X}$.
\end{itemize}
\end{thm}

The proof is based on a union of chains argument. It should be noted that for
Presburger arithmetic, $\mathrm{PR} = \Th(\omega,+)$, this argument could be
extended to any cardinality $\leq 2^{\aleph_0}$ as proved in
\cite{Knight.Nadel:82}.\footnote{For models of PR we have to define the standard
system to be the collection of sets recursive in a complete type realized in the
model, this definition is equivalent to our definition for recursively saturated
models of PA.} For PA the argument only applies for cardinalities less than or
equal $\aleph_1$, since models of PA does not admit the amalgamation property
needed, as pointed out in \cite{Knight.Nadel:82}.

For recursively saturated models the situation is almost the same as Wilmers
proved in \cite{Wilmers:75}:

\begin{thm}\label{wilmers}
For any  $\omega$-model $\scott{X}$  and complete consistent theory $T \supseteq
\PA$ the following are equivalent:
\begin{itemize}
\item[$(i)$]  $T \in \scott X$ and there is a model $M \models T$ with
$\SSy(M)=\scott{X}$, and
\item[$(ii)$] there is a recursively saturated model $M \models T$ with
$\SSy(M)=\scott{X}$.
\end{itemize}
\end{thm}

The proof is a rather straightforward application of the arithmetized
completeness theorem.

Observe that if $\scott X$ is a Scott set and $T \in \scott X$ then $\rep(T)
\subseteq \scott X$, however it may well happen that $\rep(T) \subseteq \scott
X$ and $T \notin \scott X$.

Given $(M,\scott X)$ and $I \subseteq M$ a cut, we say that $A \subseteq I$ is
coded in $\scott X$ if there is $B \in \scott X$ such that $A = B \cap
I$. Also, if $A \subseteq I \times J$ where $J$ also is a cut we say that $A$ is
coded in $\scott X$ if there is $B \in \scott X$ such that $B \cap I \times J =
A$.

Given a theory $T$ let $\Con_T$ be the first order theory consisting of the
sentences $\Con_S$, where $S\subseteq T$ ranges over all finite subtheories of
$T$ and $\Con_S$ is the sentence saying that $\lnot \sigma$ is not provable in
first order logic, where $\sigma$ is the conjunction of all sentences in $S$.

The next theorem is a special version of the arithmetized completeness theorem
tailormade for our purposes.

\begin{thm}\label{act}
Let $(M,\scott X) \models \ACA + \Con_T$, where $T$ is some first order theory
extending $\PA$ coded in $\scott X$. Then there is a model $N \models T$ with
domain $M$ and a set $C \in \scott X$ such that $N \models \varphi(a)$ iff
$\langle \varphi, a \rangle \in C$. Furthermore, the canonical embedding
$\varrho : M \to N$ is coded in $\scott X$ and such that $\varrho(M)$ is an
initial segment of $N$.
\end{thm}
\begin{proof}
Let $A \in \scott X$ code the theory $T$, i.e., $T = A \cap \omega$. By the
assumption on $(M,\scott X)$ we have $(M,\scott X) \models \Con_{x \in A \land x
< n}$ for every $n \in \omega$, where $\Con_{\theta(x)}$ is the sentence saying
that there is no proof of absurdity from sentences satisfying the formula
$\theta(x)$. 

If $M = \mathbb N$ then clearly $(M,\scott X) \models \Con_{x \in A}$ by the
compactness theorem. Assuming $M$ is nonstandard we can use overspill (this is
where we need $(M,\scott X)$ to satisfy $\ACA$) to find $T \subseteq B \in
\scott X$ such that $(M,\scott X) \models \Con_{x \in B}$.

Using the fact that the completeness theorem is provable in $\WKL$ we get the
desired set $C$ (see \cite{Simpson:99}).
\end{proof}

Let $(P,<)$ be a partial order. A filter $F\neq P$ on $P$ is an upwards
closed non-empty subset of $P$ such that if $x,y \in F$ then there is $z \in F$
satisfying $z \leq x$ and $z \leq y$. A filter $F$ is an ultrafilter if
it is a maximal filter. 

Given $(M,\scott X)$ let $P_\scott X$ be the partial order of all unbounded sets
in $\scott X$ ordered by $\subseteq$. A filter $F$ on $P_\scott X$ is
complete if for all $f : M \to M$ coded in $\scott X$ with a bounded
range there is $A \in F$ such that $f$ is constant on $A$. Any complete filter
on $P_\scott X$ is a nonprincipal ultrafilter on the boolean algebra $\scott X$.
Also, if $\scott X$ is an $\omega$-model any nonprincipal ultrafilter on $\scott
X$ is a complete filter on $P_\scott X$.

Recall that a filter $F$ on $P_\scott X$ is definable if for all $A \in \scott
X$ $\set{a  \in M  | (A)_a \in U} \in \scott X$,  where $(A)_a=\set{b \in M |
\langle a,b \rangle \in A}$ and $\langle \cdot,\cdot \rangle$ is some canonical
function coding pairs. We will, somewhat sloppily, say that a filter $U$ on
$P_\scott X$ is an ultrafilter if $U$ is an ultrafilter on $\scott X$.

\begin{lem}[{\cite{Kirby:84}}]\label{lem:def.ar}
If $(M,\scott X) \models \RCA$ and there is a definable ultrafilter $U$ on
$P_\scott X$ then $(M,\scott X) \models \ACA$. 
\end{lem}
\begin{proof}
To see this let $B = \set{a \in M| M \models \exists x \varphi(a,x,A)}$, where
$A \in \scott X$ and $\varphi$ is $\Delta_0^0$. Define $C = \set{ \langle a,b
\rangle  | \exists x \mathord< b \ \varphi (a,x,A)} \in \scott X$, then $$(C)_a
= \set{b\in M | \exists x \mathord< b \ \varphi (a,x,A)}$$ and thus $B = \set{a
\in M| (C)_a \in U} \in \scott X$ since $U$ is definable.
\end{proof}

\section{The construction}\label{sec:con}

We are now in a position to start proving the following theorem. Theorem
\ref{maintheorem} will follow from it.

\begin{thm}\label{mainlemma}
If $(M,\scott X) \models \ACA$, $P_\scott X$ has a definable complete filter $F$
and $T \supseteq \PA$ is coded in $\scott X$ such that $M \models \Con_{T}$,
then there is an end-extension $N \models T$ of $M$ satisfying $\Cod(N/M)=\scott
X$.
\end{thm}

First we construct the ultrapower we will use to build the model $N$.

Given the setup of Theorem \ref{mainlemma} let $K_0$ be the model of $T$ given
by Theorem \ref{act} and let $\varrho: K \cong K_0$ be such that $M \subseteq_e
K$ and $\varrho \upharpoonright M$ is the canonical embedding of $M$ into $K_0$.

Let $\prod_\scott{X} K$ be the set of all functions $f: M \to K$ such that the
function $\varrho \circ f: M \to M$ is coded in $\scott X$.

For any ultrafilter $U$ on $P_\scott{X}$ define $\prod_\scott{X} K / U$ to be
the set of equivalence classes of the equivalence relation $\equiv_U$ defined on
$\prod_\scott{X} K$ by $f \equiv_U g$ iff the set where $f$ and $g$ are equal,
$\set{a \in M | f(a)=g(a)}$, is in $U$. The collection $\prod_\scott{X} K / U$
of equivalence classes can be interpreted as a structure in the language of
arithmetic by the ordinary definitions of functions and relations. Let $N$
denote some model of the form $\prod_\scott{X} K /U$, where $U$ will be
understood to be an ultrafilter on $P_\scott X$. 

Let $\sigma$ be a sentence in the language $\La_A(\prod_\scott{X} K)$, i.e. the
language of arithmetic extended with the set $\prod_\scott X K$ as parameters.
The $\La_A(K)$-sentence we get by replacing all occurrences of functions $f$ by
the value $f(i)$ will be denoted by $\sigma[i]$. By $[\sigma]$ we mean the
$\La_A(N)$-sentence we get by replacing all functions $f$ by the equivalence
class $[f]$.

The \L o\'s theorem in this setting follows:

\begin{lem}
For any sentence $\sigma$ of $\La_A(\prod_\scott{X} K)$ we have $N \models
[\sigma] $ iff $$ \set{i \in M | K \models \sigma[i]} \in U.$$ 
\end{lem}

Let us embed $K$ in $N$ via the canonical embedding $F: K \to N$, $F(a) =
[f_a]$, where $f_a(b)=a$ for all $b \in M$. We will identify $a \in K$ with
$F(a) \in N$ making $K$ a substructure of $N$. In fact \L o\'s theorem gives us
that that $K \embin N$ and thus that $N \models T$. 

Let us summarize. We have $M \subseteq_e K \prec N$, where $K \models T$ and
there is a model $K_0$ with the same domain as $M$ such that $\Th(K_0,a)_{a \in
K_0}$ is coded in $\scott X$ and $\varrho : K \to K_0$ is an isomorphism. 

\begin{lem}\label{complete}
If $U$ is complete then $M \subseteq_e N$.
\end{lem}
\begin{proof}
Let $[f] \in N$ be such that $N \models [f] < a$, $a \in M$. Take $g \in [f]$
such that $g(b) < a$ for all $b \in M$.  Observe that the function $\varrho
\upharpoonright M : M \to M$ is coded in $\scott X$ and that $g(M) \subseteq M$.
Since $g(M)$ is bounded in $M$ so is $\varrho(g(M))$.

Now $\varrho \circ g : M \to M$ has bounded range and thus, by the completeness
of $U$, there is $A \in U$ such that $\varrho (g(A))=\set{b}$, $b \in M$.
Therefore $N \models [f]=[g] = \varrho^{-1}(b)$ and thus $[f] \in M$.
\end{proof}

\begin{lem}\label{lem2}
For any complete ultrafilter $U$ on $P_\scott X$, $\scott X \subseteq
\Cod(N/M)$.
\end{lem}
\begin{proof}
Given $A \in \scott{X}$ we will find $[f] \in N$ such that $[f]$ codes $A$,
i.e., $N \models a \in [f]$ iff $a \in A$, for all $a \in M$.

Let $f(a)$ be the least code in $M$ of the bounded set $A \cap I_{\mathord<a}$.
Since $M \subseteq_e K$ we have $M\prec_{\Delta_0} K$ and thus $f(a)$ also codes
the set $A \cap I_{\mathord<a}$ in $K$.

It should be clear that $f$ is coded in $\scott{X}$. For any $a,b \in M$ we have
$K \models a\in f(b) $ iff  $a < b$ and $a \in A$, and so
\begin{align*}
a \in A &\quad\Rightarrow\quad \set{ b \in I| K \models a \in f(b)} =I \setminus
 I_{\mathord<a+1}, \\
a \notin A &\quad\Rightarrow\quad \set{ b \in I | K \models a \in f(b)}  =
\emptyset.
\end{align*}

Thus, $a \in A$ iff $\set{ b \in M | K \models a \in f(b)} \in U$, i.e., iff $N
\models a \in [f]$.
\end{proof}

\begin{lem}\label{lem3}
 If $U$ is complete and definable then $\scott X \supseteq \Cod(N/M)$.
\end{lem}
\begin{proof}
We show that  $\cod_{N/M}([f]) \in \scott X$ for any $[f] \in N$. By \L o\'s's
theorem $$ \cod_{N/M}([f])  = \set{a \in M | \set{b \in M | K \models a \in
f(b)}\in U }. $$ Thus, if $A = \set{\langle a,b\rangle| K \models a \in f(b)
\text{ and } a,b \in M} \in \scott X$ then $\cod_{N/M}([f]) = \set{a \in M |
(A)_a \in U } \in \scott X$ by the definability of $U$.
\end{proof}

\begin{proof}[Proof of Theorem \ref{mainlemma}.]
Given $(M,\scott X) \models \ACA$ and $U$ a definable and complete
ultrafilter on $P_\scott X$ and a theory $T \supseteq \PA$ coded in $\scott
X$ such that $M \models \Con_{T}$, let $N$ be $\prod_\scott{X} K / U$ where $K$
is as above. By the preceding lemmas $N \models T$, $M \subseteq_e N$ and
$\Cod(N/M)=\scott X$.
\end{proof}

Now we are ready to prove the central theorem of this paper, Theorem
\ref{maintheorem}. Let us recall it.

\newcounter{tmpctr}
\setcounter{tmpctr}{\value{thm}}
\setcounter{thm}{0}
\begin{thm}
Let $\scott X$ be a Scott set that carries a definable ultrafilter and let $T
\in \scott X$ be a consistent completion of PA. Then there is a recursively
saturated model of $T$ with standard system $\scott X$.
\end{thm}
\begin{proof}
If $\scott X$ is a Scott set admitting a definable ultrafilter then $(\mathbb
N,\scott X) \models \ACA$. Any ultrafilter on $\scott X$ is complete so
Theorem \ref{mainlemma} gives us a model $N \models T$ with $\SSy(N)=\scott
X$. An application of Theorem \ref{wilmers} gives a recursively saturated model
of $T$ with standard system $\scott X$.
\end{proof}
\setcounter{thm}{\value{tmpctr}}

A natural question to ask is whatever every $\Cod(M/I)$ where $I$ is a strong
cut in $M$ admits a complete definable nonprincipal ultrafilter. This is not the
case; Enayat, in \cite{Enayat:06b}, constructed a Scott set $\scott X \models
\ACA$ of cardinality $\aleph_1$ with no nonprincipal definable ultrafilter.
However, by Theorem \ref{union} and \ref{wilmers} there is a recursively
saturated $M \models \PA$ such that $\SSy(M)=\scott X$. That $\omega$ is strong
in $M$ follows from the fact that $\scott X \models \ACA$.

\section{Definability and forcing}

Theorem \ref{maintheorem} partly reduces the question of constructing models
of
arithmetic with standard system $\scott X$ to the question when we can find
definable ultrafilters on $P_\scott X$. One way of constructing these
ultrafilters is by finding generic ultrafilters.

\begin{defin} 
Let $P$ be a partial order and $(M,\scott X) \models \WKL$.
\begin{itemize}
 \item A set $D \subseteq P$ is dense in $P$ if for every $x \in P$ there is $y
\in D$ such that $y \leq x$.
\item If $\mathscr{D}$ is a set of dense sets in $P$  a filter $F$ is
\textit{$\mathscr{D}$-generic} if $F \cap D \neq \emptyset$ for every $D \in
\mathscr{D}$.
\item A filter on $P_\scott X$ is generic if it is $\mathscr D$-generic,
where $\mathscr D$ is the collection of all dense sets (parameter) definable in
$(M,\scott X)$.
\end{itemize}
\end{defin}

It is well known that any generic filter is definable, but see \cite{Enayat:06}
for a proof.

\begin{lem}
If $(M,\scott X) \models \ACA$ and $U$ is a generic ultrafilter on $P_{\scott
X}$ then $U$ is complete and definable.
\end{lem}

A partial order $(P,<)$ is said to have the countable chain condition
(c.c.c. for short) if for every uncountable set $A \subseteq P$ there are $x,y
\in A$ such that $x$ and $y$ are compatible, i.e., there is $z \in P$ such that
$z < x$ and $z<y$.

Martin's axiom, MA for short, says that for any partial order $P$ with
the c.c.c. and any collection $\mathscr{D}$ of dense sets of cardinality $<
2^{\aleph_0}$ there is a $\mathscr{D}$-generic filter on $P$. Clearly,
$\mathrm{ZFC}+\mathrm{CH} \prf \mathrm{MA}$, but it is also the case that if
$\mathrm{ZFC}$ is consistent then so is $\mathrm{ZFC} +\mathrm{MA}+
\lnot\mathrm{CH}$. In fact, if $\mathrm{ZFC}$ is consistent and $\kappa \geq
\omega_1$ is regular such that $2^{<\kappa}=\kappa$, then 
$\mathrm{ZFC}+\mathrm{MA}+ 2^{\aleph_0}=\kappa$ is consistent.  

We can now state the following promising looking corollary:

\begin{cor}
If $\mathrm{MA}$ holds, $|\scott{X}| < 2^{\aleph_0}$ is an arithmetically
closed Scott set, i.e., $\scott X \models \ACA$, such that $P_\scott X$ has the
c.c.c. and $T \supseteq \PA$ is a consistent theory coded in $\scott X$ then
there is a recursively saturated $N \models T$ such that $\SSy(N)=\scott{X}$.
\end{cor}

However, as Hamkins and Gitman recently proved in \cite{Hamkins.Gitman:05},
there are no uncountable such sets $\scott X$.

\begin{thm}
If $(M,\scott X) \models \RCA$ is such that $P_\scott X$ has the c.c.c. then
$\scott X$ is countable.
\end{thm}
\begin{proof}
For every $A \subseteq M$, let $A^*$ be the set of codes of initial segments of
$A$, i.e., $A^* = \set{ \langle A \cap I_{< a}\rangle | a \in M}$, where
$\langle B \rangle$ is the least code in $M$ of the bounded set $B$. It should
be clear that if $X \in \scott X$ then $X^* \in \scott X$.

Also if $X \neq Y \in P_\scott{X}$ then $X^* \cap Y^*$ is bounded and thus $X^*$
and $Y^*$ are incompatible elements in $P_{\scott X}$. Therefore, the set
$\set{X^* | X \in P_\scott X }$ is an anti-chain in $P_{\scott X}$ of the same
cardinality as $\scott X$.
\end{proof}

Thus, this construction together with Martin's axiom does not give us anything
new on the existence of standard systems.

There might still be some hope of at least proving the consistency of
$\mathrm{ZFC} + \lnot \mathrm{CH} + \text{`$\ACA \subseteq \SSy$'}$, where
`$\ACA \subseteq \SSy$' denotes the sentence saying that every arithmetically
closed Scott set is the standard system of some model of PA. This might be done
by forcing generic ultrafilters into a model of $\mathrm{ZFC} + \lnot
\mathrm{CH}$. However, for this mission to succeed we have to control the
cardinals while doing the forcing. The following result of Enayat
\cite{Enayat:06b} is somewhat discouraging in that direction.

\begin{thm}
There is a Scott set $\scott X \models \ACA$ such that forcing with $P_\scott X$
collapses $\aleph_1$.
\end{thm}

Recently Gitman \cite{Gitman:06} used similar ideas to show that assuming the
proper forcing axiom (PFA) any proper arithmetically closed Scott set is a
standard system. Gitman also proved the consistency of $\lnot\mathrm{CH}$
together with the existence of an arithmetically closed proper Scott sets of
size $\aleph_1$. However it is, to our knowledge, open if PFA is consistent
with the existence of such a Scott set.

\end{document}